\documentclass[12pt]{article}
\usepackage {amsmath, amssymb, color}
\def\ep{\varepsilon}

\newtheorem{Theorem}{Theorem}[section]

\newtheorem{Remark}[Theorem]{Remark}
\newtheorem{Prop}[Theorem]{Proposition}

\begin{document}

\title 
{Spectral asymptotics for a class of integro-differential equations arising in the theory of fractional Gaussian processes}

\author{Alexander I. Nazarov\footnote{
St.Petersburg Department of Steklov Institute, Fontanka 27, St.Petersburg, 191023, Russia, 
and St.Petersburg State University, 
Universitetskii pr. 28, St.Petersburg, 198504, Russia. E-mail: al.il.nazarov@gmail.com.
}
}
\date{}

\maketitle

\begin{abstract}
\footnotesize
We study spectral problems for integro-differential equations arising in the theory of Gaussian processes similar to the fractional Brownian motion. We generalize the method of Chigansky--Kleptsyna and obtain the two-term eigenvalue asymptotics for such equations. Application to the small ball probabilities in $L_2$-norm is given.
\end{abstract}

\section{Introduction}

The spectral analysis of Gaussian processes is intensively developed in the last two decades, in particular, in the context of the problem of small deviation asymptotics for such processes in the Hilbert norm. 

It is known, see \cite{Naz09b}, that to obtain the {\it logarithmic} $L_2$-small ball asymptotics of a Gaussian process $X$, it is sufficient to know one-term asymptotics of the eigenvalue counting function of its covariance operator. However, to manage the {\it exact} asymptotics (up to a constant), we need at least two-term asymptotics of the eigenvalues with a proper remainder estimate (\cite{Li92}, 
see also \cite{GHT2003a}).

The last problem is quite delicate and was solved only for several special
processes. Most of them are the so-called {\it Green Gaussian processes}, i.e. the processes the covariance functions $G_X$ of which are the Green functions for the ordinary differential operators (ODO) subject to proper homogeneous boundary conditions. This class contains many classical processes, e.g., the integrated Brownian motion, the Slepian process and the
Ornstein-Uhlenbeck process. The special nature of the Green Gaussian processes allows us to use the well-developed techniques of spectral theory for ODOs, see, e.g., \cite{Nai}. This approach 
was elaborated in \cite{NaNi04}, \cite{Naz09a} and was used in a number of papers, see, e.g., \cite{NaNi18} and references therein. We mention also the papers \cite{Naz09}, \cite{NaPe}, \cite{Pe} 
where the two-term spectral asymptotics was obtained for finite-dimensional perturbations of the Green Gaussian processes.

The case of {\it fractional Gaussian processes} is much more complicated. Until recently,
only the main term of spectral asymptotics was known, and thus, only logarithmic small ball asymptotics was obtained for such processes. Namely, in the pioneer paper \cite{Bron} the one-term 
spectral asymptotics was calculated for the fractional Brownian motion (FBM) $W^H$, i.e. the zero mean-value Gaussian process with covariance function
\begin{equation*}
{\cal G}(x,y):=G_{W^H}(x,y)=\frac 12\,\big(x^{2H}+y^{2H}-|x-y|^{2H}\big)
\end{equation*}
(here $H\in(0,1)$ is the so-called {\it Hurst index}, the case $H=\frac 12$ corresponds to the standard Wiener process). 

A more general approach was developed in \cite{NaNi04a}. This approach is based on the powerful theorems on spectral asymptotics of integral operators \cite{BS70}, see also \cite[Appendix 7]{BS74}, 
and covers many fractional processes. However, it also gives only the one-term eigenvalue asymptotics.

A breakthrough step was managed in the paper \cite{ChiKl}. The eigenproblem for the covariance operator of $W^H$ on $[0,1]$ was reduced to the generalized eigenproblem
\begin{equation}\label{problem}
(\mathbb{K}_{\alpha}\psi)(x)=-\lambda \psi''(x),\quad x\in(0,1)
\end{equation}
with boundary conditions $\psi'(0)=\psi(1)=0$. Here we use the notation $\alpha=2-2H\in(0,2)\setminus\{1\}$, and
\begin{equation*}
(\mathbb{K}_{\alpha}\psi)(x)=(1-\alpha/2)\, \frac d{dx}\int\limits_0^1 \text{sign}(x-y)|x-y|^{1-\alpha}\psi(y)\,dy.
\end{equation*}
By the Laplace transform
\begin{equation*}
\widehat \psi(z)=\int\limits_0^1 \psi(y)\exp(-zy)\,dy
\end{equation*}
the problem (\ref{problem}) was reduced to the Riemann--Hilbert problem which, in turn, was solved asymptotically using the ideas of \cite{Uk}, \cite{Pa} (see some additional references in \cite{ChiKl}). 
In this way the two-term asymptotics of the eigenvalues with the remainder estimate was obtained for the FBM in the full range of the Hurst index. Based on this, the exact $L_2$-small ball asymptotics 
for $W^H$ was established for the first time, along with some other applications. It should be mentioned that the eigenfunction asymptotics for FBM was also obtained in \cite{ChiKl}.

In later papers \cite{ChiKlM17}, \cite{ChiKlM18} similar results were obtained for some other fractional Gaussian processes.
\medskip

In this paper we provide a slightly more general point of view and consider the eigenproblem (\ref{problem}) with general self-adjoint boundary conditions. This gives a unified approach that encompasses the previous results and covers several new fractional Gaussian processes.
\medskip

The paper is organized as follows. In Section \ref{Sec2} we calculate the two-term spectral asymptotics of the problem (\ref{problem}) with arbitrary self-adjoint boundary conditions that do not contain the spectral parameter. Here we mainly follow the line of \cite{ChiKl}. It turns out that there are three possible ``shifts'' of the second term in the asymptotics depending on the {\it sum of orders of the boundary conditions}. It is well known, see \cite[Theorem 7.1]{NaNi04} and \cite[Theorem 1.1]{Naz09a}, that this parameter drives the second term of spectral asymptotics for ODOs of arbitrary order. We conjecture that this is also the case for general eigenproblems of the type (\ref{problem}) with ODO of arbitrary order on the right-hand side.
 
 In Section \ref{Sec3} we consider a more general eigenproblem
\begin{equation}\label{problem1}
(\mathbb{K}_{\alpha}\psi)(x)=\lambda \big(-\psi''(x)+\mathfrak{p}(x)\psi(x)\big),\quad x\in(0,1),
\end{equation}
with self-adjoint boundary conditions. We prove that the additional term in (\ref{problem1}) can be considered as a weak perturbation of the problem (\ref{problem}) which does not affect the two-term 
eigenvalue asymptotics. 

In Section \ref{Sec4} we give several examples of fractional Gaussian processes covered by the results of Sections \ref{Sec2} and \ref{Sec3}.
 
Finally, in Section \ref{Sec5} we collect the results on $L_2$-small ball probabilities for the fractional processes considered in Section \ref{Sec4}.

\section{Analysis of the problem (\ref{problem}) with general boundary conditions}\label{Sec2} 

First we consider in detail the case $\alpha<1$. In this case the equation (\ref{problem}) reads as follows:
\begin{equation}\label{basic}
(1-\alpha/2)(1-\alpha) \int\limits_0^1 |x-y|^{-\alpha}\psi(y)\,dy=-\lambda \psi''(x).
\end{equation}

\subsection{Transformation of the problem}

 Following \cite[Sec. 5.1]{ChiKl} we define 
\begin{equation*}
u(x,t):=\int\limits_0^1 \exp(-t|x-y|)\psi(y)\,dy;\qquad u_0(x)=\int\limits_0^\infty t^{\alpha-1}u(x,t)\,dt.
\end{equation*}
Then (\ref{basic}) becomes
\begin{equation*}
c_\alpha u_0(x)=-\lambda\psi''(x), \qquad c_\alpha=\frac{(1-\alpha/2)(1-\alpha)}{\Gamma(\alpha)}.
\end{equation*}
The Laplace transform gives
\begin{equation*}
\widehat u_0(z)=-\frac{\lambda}{c_\alpha}\,\big(z^2\widehat\psi(z)+\exp(-z)(\psi'(1)+z\psi(1))-(\psi'(0)+z\psi(0)\big).
\end{equation*}
On the other hand,
\begin{equation*}
(z^2-t^2)\widehat u(z,t)=u(0,t)(z+t)-\exp(-z)u(1,t)(z-t)-2t\widehat\psi(z),
\end{equation*}
i.e. for $z\notin\mathbb R$
\begin{equation*}
\widehat u_0(z)=\int\limits_0^\infty \frac {t^{\alpha-1}}{z-t}\,u(0,t)\,dt-
\exp(-z)\int\limits_0^\infty \frac {t^{\alpha-1}}{z+t}\,u(1,t)\,dt-\widehat\psi(z)\int\limits_0^\infty \frac {2t^{\alpha}}{z^2-t^2}\,dt.
\end{equation*}
So, we obtain
\begin{equation*}
\aligned
\Big(\frac{\lambda}{c_\alpha}\,z^2-\int\limits_0^\infty \frac {2t^{\alpha}}{z^2-t^2}\,dt\Big)\widehat\psi(z) & =\frac{\lambda}{c_\alpha}(\psi'(0)+z\psi(0)\big)
+\int\limits_0^\infty \frac {t^{\alpha-1}}{t-z}\,u(0,t)\,dt\\
-\exp(-z)&\, \Big(\frac{\lambda}{c_\alpha}(\psi'(1)+z\psi(1)\big)-\int\limits_0^\infty \frac {t^{\alpha-1}}{z+t}\,u(1,t)\,dt\Big),
\endaligned
\end{equation*}
 and thus
\begin{equation}\label{main}
z\widehat\psi(z)=\,\frac 1{\Lambda(z)}\big(\exp(-z)\Psi(-z)+\Phi(z)\big),
\end{equation}
where
\begin{equation*}
\aligned
\Lambda(z)= &\ \frac{\lambda}{c_\alpha}\,z+\frac 1z\int\limits_0^\infty \frac {2t^{\alpha}}{t^2-z^2}\,dt\\
= &\ \frac{\lambda}{c_\alpha}\,z+z^{\alpha-2}\,\frac {\pi\exp(\pm i\pi(1-\alpha)/2)}{\cos(\pi\alpha/2)},\qquad \Im(z)\gtrless 0;
\endaligned
\end{equation*}
\begin{equation}\label{Phi}
\aligned
\Phi(z)= & \frac{\lambda}{c_\alpha}\,(\psi'(0)+z\psi(0)\big)+\int\limits_0^\infty \frac {t^{\alpha-1}}{t-z}\,u(0,t)\,dt;\\
\Psi(z)= &\ -\frac{\lambda}{c_\alpha}\,(\psi'(1)-z\psi(1)\big)+\int\limits_0^\infty \frac {t^{\alpha-1}}{t-z}\,u(1,t)\,dt.
\endaligned
\end{equation}

The function $\Lambda$ is defined in $\mathbb C\setminus\mathbb R$, has two purely imaginary zeros
\begin{equation}\label{nu}
\pm z_0=i\nu,\qquad \nu^{\alpha-3}=\frac{\lambda}{c_\alpha}\frac {\cos(\pi\alpha/2)}{\pi}
\end{equation}
and has limits on the real axis
\begin{equation*}
\Lambda^{\pm}(t):=\lim\limits_{z\to t\pm i0}\Lambda(z)=\frac{\lambda}{c_\alpha}\,t \pm |t|^{\alpha-2}
\begin{cases}
\dfrac {\pi\exp(i\pi(1\mp \alpha)/2)}{\cos(\pi\alpha/2)}, & t>0;\\
\\
\dfrac {\pi\exp(i\pi(1\pm \alpha)/2)}{\cos(\pi\alpha/2)}, & t<0.
\end{cases}
\end{equation*}
The following relations hold true:
\begin{equation}\label{limit-Lambda}
\Lambda^-(t)=\overline{\Lambda^+}(t)=-\Lambda^+(-t).
\end{equation}

We introduce the function $\theta(t)=\arg(\Lambda^+(t))=\pi-\theta(-t)$ and notice that (\ref{nu}) implies
\begin{equation}\label{theta0}
\theta_0(t):=\theta(\nu t)=\arctan\,\frac {\sin(\frac{\pi(1-\alpha)}2)}{\cos(\frac{\pi(1-\alpha)}2)+t^{3-\alpha}}, \qquad t>0.
\end{equation}
Evidently, $\theta_0$ is independent of $\nu$, positive and decreasing, $\theta_0(0+)=\frac{\pi(1-\alpha)}2$ and $\theta_0(+\infty)=0$. Moreover, integration by parts and \cite[3.252.12]{GR} give
\begin{equation*}
 \aligned
\int\limits_0^\infty \theta_0(t)\,dt= &\ \int\limits_0^\infty \frac {(\sin(\frac{\pi(1-\alpha)}2)s)^{\frac 1{3-\alpha}}}{s^2+2s\cot(\frac{\pi(1-\alpha)}2)+\csc^2(\frac{\pi(1-\alpha)}2)}\\
=&\ \pi\,\frac {\sin(\frac{\pi(1-\alpha)}{2(3-\alpha)})}{\sin(\frac{\pi}{3-\alpha})}=\pi\cot\big(\frac{\pi}{3-\alpha}\big)=:\pi\,\mathfrak{b}_{\alpha}.
 \endaligned
\end{equation*}

Now we look at the equation (\ref{main}) on the real line. It shows that the right-hand side is continuous on $\mathbb{R}$, and we obtain for $t>0$ and $t<0$ respectively
\begin{equation*}
 \aligned
 \frac 1{\Lambda^+(t)}\big(\exp(-t)\Psi(-t)+\Phi^+(t)\big)= \frac 1{\Lambda^-(t)}\big(\exp(-t)\Psi(-t)+\Phi^-(t)\big);\\
  \frac 1{\Lambda^+(t)}\big(\exp(-t)\Psi^-(-t)+\Phi(t)\big)= \frac 1{\Lambda^-(t)}\big(\exp(-t)\Psi^+(-t)+\Phi(t)\big),
 \endaligned
\end{equation*}
or, equivalently, with regard to (\ref{limit-Lambda}),
\begin{equation}\label{contin}
 \begin{cases}
 \Phi^+(t)-\dfrac {\Lambda^+(t)}{\Lambda^-(t)}\Phi^-(t)=\exp(-t)\Psi(-t)\Big(\dfrac {\Lambda^+(t)}{\Lambda^-(t)}-1\Big);\\
 \\
 \Psi^+(t)-\dfrac {\Lambda^+(t)}{\Lambda^-(t)}\Psi^-(t)=\exp(-t)\Phi(-t)\Big(\dfrac {\Lambda^+(t)}{\Lambda^-(t)}-1\Big), 
  \end{cases}
  \qquad t>0.
 \end{equation}
Since 
\begin{equation*}
\frac {\Lambda^+(t)}{\Lambda^-(t)}=\frac {\Lambda^+(t)}{\overline{\Lambda^+}(t)}=\exp(2i\theta(t)),
 \end{equation*}
we can rewrite (\ref{contin}) as follows:
\begin{equation}\label{contin1}
 \begin{cases}
 \Phi^+(t)-\exp(2i\theta(t))\Phi^-(t)=2i\exp(-t)\exp(i\theta(t))\sin(\theta(t))\Psi(-t);\\
 \\
 \Psi^+(t)-\exp(2i\theta(t))\Psi^-(t)=2i\exp(-t)\exp(i\theta(t))\sin(\theta(t))\Phi(-t). 
  \end{cases}
\end{equation}

We also know from definition that $\Phi(z)$ and $\Psi(z)$ behave as polynomials of order not greater than one at infinity whereas they are $O(z^{\alpha-1})$ at the origin.

We introduce the function $X_0(z)$ with the cut at the positive semiaxis such that
\begin{equation}\label{X0}
 \frac {X_0^+(t)}{X_0^-(t)}=\exp(2i\theta_0(t)),\quad t>0;\qquad X_0(z)\asymp 
 \begin{cases}
 1,& z\to\infty;\\
z^{\frac {\alpha-1}2},& z\to0.
 \end{cases}
\end{equation}
The first relation in (\ref{X0}) is satisfied by the Sokhotski--Plemelj formula
\begin{equation}\label{X0a}
 X_0(z):=\exp\Big(\frac 1{\pi} \int\limits_0^\infty\frac {\theta_0(s)}{s-z}ds\,\Big).
\end{equation}
It is easy to see that 
\begin{equation}\label{X0b}
 \aligned
 X_0(z)= &\ \exp\Big(-\frac {\mathfrak{b}_{\alpha}}z+O\big(\frac 1{z^2}\big)\Big) =1-\frac {\mathfrak{b}_{\alpha}}z+O\big(\frac 1{z^2}\big), & z\to\infty;\\
 X_0(z)\asymp &\ \exp\Big(-\frac {\theta_0(0+)}{\pi}\log(z)\Big)=z^{\frac {\alpha-1}2}, & z\to0.
 \endaligned
\end{equation}

Using (\ref{X0}) we rewrite (\ref{contin1}) as follows:
\begin{equation}\label{Phi0}
\begin{cases}
 \dfrac {\Phi_0^+(t)}{X_0^+(t)}-\dfrac {\Phi_0^-(t)}{X_0^-(t)}=2i\exp(-\nu t)\exp(i\theta_0(t))\sin(\theta_0(t))\dfrac {X_0(-t)}{X_0^+(t)}\dfrac {\Psi_0(-t)}{X_0(-t)};\\
 \\
 \dfrac {\Psi_0^+(t)}{X_0^+(t)}-\dfrac {\Psi_0^-(t)}{X_0^-(t)}=2i\exp(-\nu t)\exp(i\theta_0(t))\sin(\theta_0(t))\dfrac {X_0(-t)}{X_0^+(t)}\dfrac {\Phi_0(-t)}{X_0(-t)},
  \end{cases}
\end{equation}
where $\Phi_0(t)=\Phi(\nu t)$ and $\Psi_0(t)=\Psi(\nu t)$. Therefore, functions
\begin{equation*}
 S(z)=\frac {\Phi_0(z)+\Psi_0(z)}{2X_0(z)},\qquad  D(z)=\frac {\Phi_0(z)-\Psi_0(z)}{2X_0(z)}
\end{equation*}
satisfy for $t>0$ the equations
\begin{equation*}
 \aligned
 S^+(t)-S^-(t) & =2i\exp(-\nu t)h_0(t)S(-t);\\
  D^+(t)-D^-(t) & =-2i\exp(-\nu t)h_0(t)D(-t),
 \endaligned
\end{equation*}
where
\begin{equation*}
 \aligned
 h_0(t)= &\ \exp(i\theta_0(t))\sin(\theta_0(t))\frac {X_0(-t)}{X_0^+(t)}\\
 = &\ \sin(\theta_0(t))\exp\Big(-\frac 1{\pi} \int\limits_0^\infty\theta'_0(s)\log\Big|\frac {s+t}{s-t}\Big|ds\Big)
 \endaligned
\end{equation*}
(here we used (\ref{X0a}) and integration by parts).

Since $S(z)$ and $D(z)$ behave as polynomials of order not greater than one at infinity, the Sokhotski--Plemelj formula yields
\begin{equation*}
\aligned
 S(z) & =\frac 1{\pi} \int\limits_0^\infty\frac {\exp(-\nu s)h_0(s)}{s-z}\,S(-s)\,ds+C_1+C_2z;\\
  D(z) & =-\frac 1{\pi} \int\limits_0^\infty\frac {\exp(-\nu s)h_0(s)}{s-z}\,D(-s)\,ds+C_3+C_4z.
\endaligned
\end{equation*}

Substituting $z=-t$, $t>0$, we obtain the integral equations
\begin{equation*}
 \widehat S(t)-({\cal A} \widehat S)(t)=C_1-C_2t;\qquad
  \widehat D(t)+({\cal A} \widehat D)(t)=C_3-C_4t,
\end{equation*}
where $\widehat S(t)=S(-t)$, $\widehat D(t)=D(-t)$, and ${\cal A}$ is the integral operator with the kernel $A(t,s)=\frac {\exp(-\nu s)h_0(s)}{\pi(s+t)}$, $s,t\in \mathbb{R}_+$.

By \cite[Lemma 5.6]{ChiKl}, if $\nu$ is large enough then the operator ${\cal A}$ is contracting on $L_2(\mathbb{R}_+)$ and maps any polynomial into $L_2(\mathbb{R}_+)$. Therefore, these equations 
are uniquely solvable. Moreover, the relation $h_0(0)=\sin(\theta_0(0+))=\sin(\frac{\pi(1-\alpha)}2)$ shows that (see \cite[3.241.2]{GR})
\begin{equation*}
\widehat S(t), \widehat D(t)=O(t^{\frac{\alpha-1}2}) \ \  \mbox{as} \ \ t\to0 \quad
\Longrightarrow \quad S(z), D(z)=O(z^{\frac{\alpha-1}2}) \ \  \mbox{as} \ \ z\to0,
\end{equation*}
and therefore, (\ref{X0b}) implies $\Phi_0(z),\Psi_0(z)=O(z^{\alpha-1})$ as $z\to0$, as required.

Denote by $p_{\pm}^0(t)$ and $p_{\pm}^1(t)$ the (unique) solutions of the equations on $\mathbb R_+$
\begin{equation*}
 p_{\pm}^0(t)\mp({\cal A} p_{\pm}^0)(t)=1;\qquad p_{\pm}^1(t)\mp({\cal A} p_{\pm}^1)(t)=t,
\end{equation*}
and extend them analytically to $\mathbb{C}\setminus\mathbb{R}_-$. Then evidently
\begin{equation*}
 S(z)=C_1p_+^0(-z)-C_2p_+^1(-z),\qquad D(z)=C_3p_-^0(-z)-C_4p_-^1(-z),
\end{equation*}
whence
\begin{equation}\label{infty}
\aligned
\Phi_0(z)= &\ X_0(z)(C_1p_+^0(-z)-C_2p_+^1(-z)+C_3p_-^0(-z)-C_4p_-^1(-z));\\
 \Psi_0(z)= &\ X_0(z)(C_1p_+^0(-z)-C_2p_+^1(-z)-C_3p_-^0(-z)+C_4p_-^1(-z)).
\endaligned
\end{equation}

Since $\widehat\psi$ is an entire function, the following relation is fulfilled:
\begin{equation}\label{eigen-equ}
\exp(-z_0)\Psi(-z_0)+\Phi(z_0)\equiv\exp(-i\nu)\Psi_0(-i)+\Phi_0(i)=0, 
\end{equation}
where $z_0$ is introduced in (\ref{nu}). So, every eigenvalue of the original problem generates a root of (\ref{eigen-equ}) by the relation (\ref{nu}). By following the corresponding argument in \cite[Lemma 5.3]{ChiKl} one finds, vice versa, that every root of (\ref{eigen-equ}) (except for $\nu=0$, if it arises) generates an eigenvalue through (\ref{nu}).

We multiply (\ref{eigen-equ}) by $\exp(i\nu/2)$ and obtain
\begin{equation}\label{eigen}
\aligned
 & C_1(\exp(i\nu/2)X_0(i)p_+^0(-i)+\exp(-i\nu/2)X_0(-i)p_+^0(i))\\
 - &\ C_2(\exp(i\nu/2)X_0(i)p_+^1(-i)+\exp(-i\nu/2)X_0(-i)p_+^1(i))\\
 + &\ C_3(\exp(i\nu/2)X_0(i)p_-^0(-i)-\exp(-i\nu/2)X_0(-i)p_-^0(i))\\
 - &\ C_4(\exp(i\nu/2)X_0(i)p_-^1(-i)-\exp(-i\nu/2)X_0(-i)p_-^1(i)) =0.
\endaligned
\end{equation}
By \cite[Lemma 5.5]{ChiKl} we have 
\begin{equation*}
X_0(\pm i)=\sqrt{\frac{3-\alpha}2}\exp(\pm i\pi(1-\alpha)/8),
\end{equation*}
and \cite[Lemma 5.7]{ChiKl} claims
\begin{equation*}
\aligned
& p_{\pm}^0(i)=1+O(\nu^{-1}),\quad p_{\pm}^0(-i)=1+O(\nu^{-1}), \\
& p_{\pm}^1(i)=i+O(\nu^{-2}),\quad p_{\pm}^1(-i)=-i+O(\nu^{-2}),
\endaligned
\qquad \mbox{as}\ \ \nu\to\infty.
\end{equation*}
Thus, (\ref{eigen}) is equivalent to
\begin{equation*}
\aligned
& C_1\Big[\cos\Big(\frac {\nu+\rho}2\Big)\Big]-C_2\Big[\sin\Big(\frac {\nu+\rho}2\Big)\Big]\\
+i\ & \Big(C_3\Big[\sin\Big(\frac {\nu+\rho}2\Big)\Big]+C_4\Big[\cos\Big(\frac {\nu+\rho}2\Big)\Big]\Big)
=0
\endaligned
\end{equation*}
(here and elsewhere $\rho=\pi(1-\alpha)/4$ and we use the notation $[a]=a+O(\nu^{-1})$, see  \cite[\S4]{Nai}).

By (\ref{Phi}) and (\ref{X0b}), all coefficients $C_j$ are real, therefore, (\ref{eigen}) is equivalent to the real system
\begin{equation}\label{eigen2}
\aligned
& C_1\Big[\cos\Big(\frac {\nu+\rho}2\Big)\Big]-C_2\Big[\sin\Big(\frac {\nu+\rho}2\Big)\Big]=0;\\
& C_3\Big[\sin\Big(\frac {\nu+\rho}2\Big)\Big]+C_4\Big[\cos\Big(\frac {\nu+\rho}2\Big)\Big]=0.
\endaligned
\end{equation}

Now we compare (\ref{Phi}) and the behavior of $\Phi(\nu z)$ and $\Psi(\nu z)$ at infinity provided by (\ref{infty}). By \cite[Lemma 5.7]{ChiKl} we have
\begin{equation*}
p_{\pm}^0(z)=1+O(z^{-1}),\qquad p_{\pm}^1(z)=z+O(z^{-1}),\qquad \mbox{as}\ \ z\to\infty.
\end{equation*}
Using (\ref{Phi}) and (\ref{X0b}) we obtain
\begin{equation*}
\aligned
C_1-C_2\mathfrak{b}_{\alpha}+C_3-C_4\mathfrak{b}_{\alpha}=\frac{\lambda}{c_\alpha}\,\psi'(0); \qquad C_2+C_4=\frac{\lambda}{c_\alpha}\,\psi(0)\nu;\\
C_1-C_2\mathfrak{b}_{\alpha}-C_3+C_4\mathfrak{b}_{\alpha}=-\frac{\lambda}{c_\alpha}\,\psi'(1); \qquad C_2-C_4=\frac{\lambda}{c_\alpha}\,\psi(1)\nu.
\endaligned
\end{equation*}
We solve these equations and substitute $C_j$ into (\ref{eigen2}). This gives
\begin{equation}\label{eigen3}
\aligned
& (\psi'(0)-\psi'(1)) \big[\mathfrak{A}\big]+(\psi(0)+\psi(1))\nu\big(\mathfrak{b}_{\alpha}\big[\mathfrak{A}\big]-\big[\mathfrak{B}\big]\big)=0;\\
& (\psi'(0)+\psi'(1)) \big[\mathfrak{B}\big]+(\psi(0)-\psi(1))\nu\big(\big[\mathfrak{A}\big]+\mathfrak{b}_{\alpha}\big[\mathfrak{B}\big]\big)=0
\endaligned
\end{equation}
(here $\mathfrak{A}=\cos\big(\frac {\nu+\rho}2\big)$ and $\mathfrak{B}=\sin\big(\frac {\nu+\rho}2\big)$).

The equations (\ref{eigen3}) complemented by the boundary conditions of the original problem generate a ($4\times4$) homogeneous system. Standard argument based on the Rouch\'{e} theorem shows that 
the roots of its determinant are approximations of the solutions of (\ref{eigen-equ}) for large $|\nu|$.

\subsection{Separated boundary conditions}

Separated boundary conditions (or Sturm type conditions) for the second order operator can be written as follows:\footnote{Recall that here we consider the boundary conditions that do not contain 
the spectral parameter $\lambda$.}
\begin{equation}\label{separated}
\beta_0\psi'(0)-\gamma_0\psi(0)=0;\qquad \beta_1\psi'(1)+\gamma_1\psi(1)=0
\end{equation}
(one of two coefficients in every condition may vanish).

We denote by $\varkappa$ the sum of orders of the derivatives in boundary conditions (\ref{separated}). It was mentioned in the Introduction that this quantity plays an important role in the spectral 
asymptotics of ordinary differential operators, see, e.g., \cite{NaNi04}. In our case, evidently, 
$\varkappa\in\{0,1,2\}$.
\medskip

{\bf 1}. Let $\varkappa=0$. Then (\ref{separated}) reads $\psi(0)=\psi(1)=0$, and (\ref{eigen3}) is reduced to
\begin{equation*}
\big[\mathfrak{A}\big]\psi'(0)-\big[\mathfrak{A}\big]\psi'(1) =0;\qquad
\big[\mathfrak{B}\big]\psi'(0)+\big[\mathfrak{B}\big]\psi'(1)=0.
\end{equation*}

The existence condition for a nontrivial solution to this system takes the form
\begin{equation*}
2\mathfrak{A}\mathfrak{B}\equiv \sin(\nu+\rho)=O(\nu^{-1}),\qquad \mbox{as}\ \ \nu\to\infty.
\end{equation*}
Thus, if we enumerate the roots of (\ref{eigen-equ}) in increasing order of absolute values then 
\begin{equation}\label{kappa=0}
\nu_{n+k}=\pi n-\frac {\pi(1-\alpha)}4+O(n^{-1}),\qquad \mbox{as}\ \ n\to\infty
\end{equation}
for some $k$. 

By verbatim repetition of \cite[Section 5.1.7]{ChiKl} we show that $k$ is independent of $\alpha$. Therefore, it can be calculated by considering the case $\alpha=1$ where the original problem 
becomes the standard Sturm--Liouville problem. Thus, we obtain $k=0$.
\medskip

{\bf 2}. Let $\varkappa=1$. By symmetry we can suppose without loss of generality that (\ref{separated}) reads $\psi(0)=\psi'(1)+\gamma\psi(1)=0$, and (\ref{eigen3}) is reduced to
\begin{equation*}
\aligned
& \big[\mathfrak{A}\big]\psi'(0)+\nu\big((\mathfrak{b}_{\alpha}+\gamma\nu^{-1})\big[\mathfrak{A}\big]-\big[\mathfrak{B}\big]\big)\psi(1) =0;\\
& \big[\mathfrak{B}\big]\psi'(0)-\nu\big(\big[\mathfrak{A}\big]+(\mathfrak{b}_{\alpha}+\gamma\nu^{-1})\big[\mathfrak{B}\big]\big)\psi(1) =0.
\endaligned
\end{equation*}
The existence condition for a nontrivial solution to this system takes the form
\begin{equation*}
\mathfrak{A}^2-\mathfrak{B}^2 +2\mathfrak{b}_{\alpha}\mathfrak{A}\mathfrak{B}\equiv \cos(\nu+\rho)+\mathfrak{b}_{\alpha}\sin(\nu+\rho)=O(\nu^{-1}).
\end{equation*}
Recalling that $\mathfrak{b}_{\alpha}=\cot\big(\frac{\pi}{3-\alpha}\big)$ we conclude that
in this subcase
\begin{equation}\label{kappa=1}
\nu_{n+k}=\pi n-\frac {\pi(1-\alpha)}4-\frac {\pi}{3-\alpha}+O(n^{-1}),\qquad \mbox{as}\ \ n\to\infty
\end{equation}
for some $k$. Comparing with the case $\alpha=1$ we obtain $k=0$.
\medskip

{\bf 3}. Let $\varkappa=2$. In this case (\ref{separated}) reads $\psi'(0)-\gamma_0\psi(0)=\psi'(1)+\gamma_1\psi(1)=0$, and (\ref{eigen3}) is reduced to
\begin{equation*}
\aligned
& \big((\mathfrak{b}_{\alpha}+\gamma_0\nu^{-1})\big[\mathfrak{A}\big]-\big[\mathfrak{B}\big]\big)\psi(0)+\big((\mathfrak{b}_{\alpha}+\gamma_1\nu^{-1})\big[\mathfrak{A}\big]-\big[\mathfrak{B}\big]\big)\psi(1) =0;\\
& \big(\big[\mathfrak{A}\big]+(\mathfrak{b}_{\alpha}+\gamma_0\nu^{-1})\big[\mathfrak{B}\big]\big)\psi(0)-\big(\big[\mathfrak{A}\big]+(\mathfrak{b}_{\alpha}+\gamma_1\nu^{-1})\big[\mathfrak{B}\big]\big)\psi(1) =0.
\endaligned
\end{equation*}
The existence condition for a nontrivial solution to this system takes the form
\begin{equation*}
\mathfrak{b}_{\alpha}\big(\mathfrak{A}^2-\mathfrak{B}^2\big) +(\mathfrak{b}_{\alpha}^2-1)\mathfrak{A}\mathfrak{B}\equiv \mathfrak{b}_{\alpha}\big(\cos(\nu+\rho)
+\cot\big(\frac{2\pi}{3-\alpha}\big)\sin(\nu+\rho)\big)=O(\nu^{-1}),
\end{equation*}
and we conclude that in this subcase
\begin{equation}\label{kappa=2}
\nu_{n+k}=\pi n-\frac {\pi(1-\alpha)}4-\frac {2\pi}{3-\alpha}+O(n^{-1}),\qquad \mbox{as}\ \ n\to\infty
\end{equation}
for some $k$. Comparing with the case $\alpha=1$ we obtain $k=0$. 
\medskip

Now we can formulate the final result.

\begin{Theorem}
 The eigenvalues of the problem (\ref{basic}) with separated boundary conditions have the following asymptotics as $n\to\infty$:
\begin{equation}\label{separated-final}
\lambda_n=\sin(\pi\alpha/2)\Gamma(3-\alpha) \Big(\pi n-\frac {\pi(1-\alpha)}4-\frac {\varkappa\pi}{3-\alpha}+O(n^{-1})\Big)^{\alpha-3},
\end{equation}
where $\varkappa$ stands for the sum of orders of the derivatives in conditions (\ref{separated}).
\end{Theorem}

This statement easily follows from the relation (\ref{nu}) and the obtained asymptotics of $\nu_n$.

\begin{Remark}\label{remark-separ}
Notice that in general zero root of (\ref{eigen-equ}) may arise (say, for $\varkappa=2$, $\gamma_0=\gamma_1=0$). Since $\nu=0$ does not generate an eigenvalue by formula (\ref{nu}), 
this forces us to shift by one the numbering in (\ref{separated-final}).
\end{Remark}

\subsection{Almost separated boundary conditions}

For the second order operator, almost separated (or separated in the principal order) boundary conditions can be written as follows:
\begin{equation*}
\psi'(0)-\gamma_0\psi(0)-\widehat\gamma\psi(1)=0;\qquad \psi'(1)+\gamma_1\psi(1)+\widehat\gamma\psi(0)=0.
\end{equation*}
Analysis of this case mostly repeats the subcase $\varkappa=2$, and the eigenvalue asymptotics coincides with (\ref{separated-final}) for $\varkappa=2$.

\begin{Remark}
In this case zero root of (\ref{eigen-equ}) may arise (say, for $\gamma_0=\gamma_1=-\widehat\gamma$), which forces us to shift the numbering by one.
\end{Remark}

\subsection{Non-separated boundary conditions}

For the second order operator, non-separated boundary conditions can be written as follows:
\begin{equation}\label{non-separated}
\beta\psi'(0)+\gamma\psi'(1)+\delta\psi(0)=0;\qquad \gamma\psi(0)+\beta\psi(1)=0.
\end{equation}
In this case (\ref{eigen3}) is reduced to
\begin{equation*}
\aligned
& \big(1+\frac {\beta}{\gamma}\big)\big[\mathfrak{A}\big]\psi'(0)+\nu\Big(\big(1-\frac {\gamma}{\beta}\big)\big(\mathfrak{b}_{\alpha}\big[\mathfrak{A}\big]-\big[\mathfrak{B}\big]\big)
+\frac {\delta}{\gamma}\,\nu^{-1}\big[\mathfrak{A}\big]\Big)\psi(0) =0;\\
& \big(1-\frac {\beta}{\gamma}\big)\big[\mathfrak{B}\big]\psi'(0)+\nu\Big(\big(1+\frac {\gamma}{\beta}\big)\big(\big[\mathfrak{A}\big]+\mathfrak{b}_{\alpha}\big[\mathfrak{B}\big]\big)
-\frac {\delta}{\gamma}\,\nu^{-1}\big[\mathfrak{B}\big]\Big)\psi(0) =0.
\endaligned
\end{equation*}
The existence condition for a nontrivial solution to this system takes the form
\begin{equation*}
\aligned
& \frac {(\beta+\gamma)^2}{\beta^2+\gamma^2}\,\mathfrak{A}^2-\frac {(\beta-\gamma)^2}{\beta^2+\gamma^2}\,\mathfrak{B}^2 +2\mathfrak{b}_{\alpha}\mathfrak{A}\mathfrak{B}\\
\equiv\ & \cos(\nu+\rho)+\mathfrak{b}_{\alpha}\sin(\nu+\rho)+\frac {2\beta\gamma}{\beta^2+\gamma^2}=O(\nu^{-1}).
\endaligned
\end{equation*}
Therefore, in this case the sequence $\nu_n$ can be split into two subsequences $\nu'_n$, $\nu''_n$ such that, as $n\to\infty$,
\begin{equation}\label{nu-non-separ}
\aligned
\nu'_{n+k'}= &\ \pi(2n-1)-\frac {\pi(1-\alpha)}4-\frac {\pi}{3-\alpha}\\
+ &\ \arcsin\Big(\frac {2\beta\gamma}{\beta^2+\gamma^2}\sin\Big(\frac {\pi}{3-\alpha}\Big)\Big)+O(n^{-1});\\
\nu''_{n+k''}= &\ 2\pi n-\frac {\pi(1-\alpha)}4-\frac {\pi}{3-\alpha}\\
- &\ \arcsin\Big(\frac {2\beta\gamma}{\beta^2+\gamma^2}\sin\Big(\frac {\pi}{3-\alpha}\Big)\Big)+O(n^{-1}),
\endaligned
\end{equation}
for some $k'$, $k''$. Comparing with the case $\alpha=1$ we obtain $k'+k''=0$, so without loss of generality we can put $k'=k''=0$.
\medskip

Now we can formulate the final result.

\begin{Theorem}
 The eigenvalues of the problem (\ref{basic}) with non-separated boundary conditions have the following asymptotics as $n\to\infty$:
\begin{equation}\label{non-separ-final}
\aligned
\lambda_n= &\ \sin(\pi\alpha/2)\Gamma(3-\alpha) 
\Big(\pi n-\frac {\pi(1-\alpha)}4-\frac {\pi}{3-\alpha}\\
- &\ (-1)^n\arcsin\Big(\frac {2\beta\gamma}{\beta^2+\gamma^2}\sin\Big(\frac {\pi}{3-\alpha}\Big)\Big)+O(n^{-1})\Big)^{\alpha-3}.
\endaligned
\end{equation}
\end{Theorem}

This statement easily follows from relations (\ref{nu}) and (\ref{nu-non-separ}).

\begin{Remark}\label{remark-non-separ}
In the case $\beta\gamma=0$ the boundary conditions (\ref{non-separated}) are in fact separated, with $\varkappa=1$. So, two subsequences (\ref{nu-non-separ}) can be merged, and the result 
coincides with (\ref{kappa=1}). In general case two subsequences of eigenvalues have the opposite shifts with respect to (\ref{kappa=1}), cf. \cite[Theorem 1.1]{Naz09a}.

In the cases $\beta=\pm\gamma$ one of the subsequences (\ref{nu-non-separ}) has the second term as in (\ref{kappa=0}) while the other one has the second term as in (\ref{kappa=2}). Notice that 
in contrast to the case $\alpha=1$ two subsequences in (\ref{nu-non-separ}) cannot be asymptotically close or coincide.

Also a zero root of (\ref{eigen-equ}) may arise (say, for $\beta=-\gamma$, $\delta=0$), which forces us to shift the numbering by one.
\end{Remark}

\subsection{The case $\alpha>1$}

Repeating the argument of \cite[Sec. 5.2]{ChiKl}, we arrive at the relation (\ref{main}) with
\begin{equation*}
\aligned
\Lambda(z)= &\ \frac{\lambda}{|c_\alpha|}\,z-z\int\limits_0^\infty \frac {2t^{\alpha-2}}{t^2-z^2}\,dt\\
= &\ \frac{\lambda}{|c_\alpha|}\,z+z^{\alpha-2}\,\frac {\pi\exp(\pm i\pi(1-\alpha)/2)}{|\cos(\pi\alpha/2)|},\qquad \Im(z)\gtrless 0;
\endaligned
\end{equation*}
\begin{equation*}
\aligned
\Phi(z)= & \frac{\lambda}{|c_\alpha|}\,(\psi'(0)+z\psi(0)\big)+\int\limits_0^\infty \frac {t^{\alpha-1}}{t-z}\,u(0,t)\,dt;\\
\Psi(z)= &\ -\frac{\lambda}{|c_\alpha|}\,(\psi'(1)-z\psi(1)\big)-\int\limits_0^\infty \frac {t^{\alpha-1}}{t-z}\,u(1,t)\,dt.
\endaligned
\end{equation*}
Following the same line as in previous subsections we again obtain formulae (\ref{separated-final}) and (\ref{non-separ-final}).

\section{A more general problem}\label{Sec3}

As it was explained in the Introduction, we wish to consider the problem (\ref{problem1}) as a perturbation of the problem (\ref{problem}).  For simplicity only, we assume that the operator $-\psi''$ 
with given boundary conditions is positive definite, otherwise the argument should be changed in a standard way.

We begin with the estimate for the eigenfunctions of the problem (\ref{problem}) with arbitrary self-adjoint boundary conditions that do not contain the spectral parameter. 

As in Section \ref{Sec2}, we consider the case $\alpha<1$; for $\alpha>1$ the argument is similar, and the result is the same.
Following the proof of \cite[Sec. 5.1.5]{ChiKl} we write
\begin{equation*}
\widehat{\psi'}(z)= z\widehat\psi(z)-\psi(0)+\exp(-z)\psi(1).
\end{equation*}
By construction, $\widehat{\psi'}(z)$ is an entire function, and we can restore $\psi'$ by integrating over the imaginary axis. Using the relation (\ref{main}) we obtain
\begin{equation}\label{invLapl}
\psi'(x)= -\frac 1{2\pi i}\lim\limits_{R\to\infty}\int\limits_{-iR}^{iR}\big({\cal F}_1(z)\exp(z(x-1))+{\cal F}_2(z)\exp(zx)\big)\,dz,
\end{equation}
where
\begin{equation*}
{\cal F}_1(z)=\frac {\Psi(-z)}{\Lambda(z)}+\psi(1),\qquad {\cal F}_2(z)=\frac {\Phi(z)}{\Lambda(z)}-\psi(0).
\end{equation*}
The integral in (\ref{invLapl}) does not depend on the constant terms in ${\cal F}_1$ and ${\cal F}_2$. It was calculated in the proof of \cite[Lemma 5.8]{ChiKl}. Up to a multiplicative constant, we obtain
\begin{equation}\label{eigenfun}
\aligned
&\psi'(x)=-\frac 2{3-\alpha}\,\Re\big(\exp(i\nu x)\Phi_0(i)\big)\\
&+\frac 1{\pi}\int\limits_0^{\infty} \frac {\sin(\theta_0(t))}{\tau_0(t)}\, \Big(\exp(-\nu t(1-x))\Psi_0(-t)-\exp(-\nu tx)\Phi_0(-t)\Big)dt,
\endaligned
\end{equation}
where $\nu$ is related to the eigenvalue $\lambda$ by (\ref{nu}), $\theta_0$ is introduced in (\ref{theta0}), $\Phi_0$ and $\Psi_0$ are defined after (\ref{Phi0}), and
\begin{equation*}
 \tau_0(t)=\frac {\cos(\pi\alpha/2)}{\pi\nu^{\alpha-2}}\,|\Lambda_+(\nu t)|=|t+t^{\alpha-2}\exp(i\pi(1-\alpha)/2)|.
\end{equation*}

Taking into account the behavior of the terms in (\ref{eigenfun}) at zero and infinity we derive\footnote{We limit ourselves to the eigenfunction estimate though its asymptotics can be also obtained 
from (\ref{eigenfun}) as it is done in \cite{ChiKl}--\cite{ChiKlM18}. In particular, the phase shift $\phi$ can be written explicitly. 
}
\begin{equation}\label{eigenfun1}
\psi'(x)= A\big(\cos(\nu x+\phi(\nu,\alpha))+F(\nu,\alpha,x)\big),
\end{equation}
where
\begin{equation*}
|F(\nu,\alpha,x)|\le \int\limits_0^{\infty} \frac {c(\alpha)t}{1+t^{4-\alpha}}\,(\exp(-\nu t(1-x))+\exp(-\nu tx))\,dt.
\end{equation*}
It is easy to see that
\begin{equation}\label{pogransloy}
\int\limits_0^1|F(\nu,\alpha,x)|\,dx\le c_1(\alpha)\nu^{-1}.
\end{equation}

Now we introduce the energy space ${\cal H}$ of the problem (\ref{problem}). For smooth functions $h_1$ and $h_2$ satisfying boundary conditions, we set
\begin{equation*}
(h_1,h_2)_{\cal H}:=-\int\limits_0^1h''_1(x)h_2(x)\,dx
\end{equation*}
and define the Hilbert space ${\cal H}$ as the completion of the set of such functions with respect to the norm generated by the scalar product $(\cdot,\cdot)_{\cal H}$. It is well known that, 
depending on the boundary conditions, ${\cal H}$ coincides either with standard Sobolev space $W^1_2(0,1)$ or with a codimension $1$ or $2$ subspace thereof. Corresponding norm is given by
\begin{equation}\label{norma}
\|h\|_{\cal H}^2=\int\limits_0^1(h'(x))^2\,dx+q(h,h),
\end{equation}
where $q(h,h)$ is a quadratic form of the variables $h(0)$ and $h(1)$.
\medskip

In a standard way, we rewrite the problems (\ref{problem}) and (\ref{problem1}) as the equations in ${\cal H}$
\begin{equation}\label{equations}
{\cal K}\psi=\lambda \psi;\qquad {\cal K}\psi=\lambda (\psi+{\cal B}\psi),
\end{equation}
where ${\cal K}$ and ${\cal B}$ are compact self-adjoint operators in ${\cal H}$ defined by the relations
\begin{equation*}
({\cal K}\psi,\eta)_{\cal H}:=\int\limits_0^1(\mathbb{K}_{\alpha}\psi)(x)\eta(x)\,dx;\quad ({\cal B}\psi,\eta)_{\cal H}:=\int\limits_0^1\mathfrak{p}(x)\psi(x)\eta(x)\,dx.
\end{equation*}

We formulate an important abstract statement.

\begin{Prop}\label{prop} (Theorem 1 in \cite{Naz19}). 
 Let ${\cal K}$ and ${\cal B}$ be self-adjoint compact operators in the Hilbert space ${\cal H}$. Suppose that ${\cal K}$ and ${\cal I}+{\cal B}$ are positive. Denote by $\lambda_n$ the eigenvalues 
 of ${\cal K}$ enumerated in the decreasing order taking into account the multiplicities, and by $\psi_n$ corresponding normalized eigenfunctions. Finally, suppose that
\begin{equation}\label{asymp}
\lambda_n=\big(an+b+O(n^{-\ep})\big)^{-r}, \qquad  \mbox{as}\quad n\to\infty,
\end{equation}
\begin{equation}\label{perturb}
 |({\cal B}\psi_n,\psi_m)_{\cal H}|\le c(mn)^{-\frac {1+\ep}2},
\end{equation}
where $a,c,\ep,r>0$, $b\in\mathbb{R}$. 
 Then the eigenvalues $\mbox{\boldmath$\lambda$}_n$ of the generalized eigenproblem
\begin{equation*}
{\cal K}\mbox{\boldmath$\psi$}_n=\mbox{\boldmath$\lambda$}_n\big(\mbox{\boldmath$\psi$}_n+{\cal B}\mbox{\boldmath$\psi$}_n\big)
\end{equation*}
 have the same two-term asymptotics as $n\to\infty$:
 \begin{equation*}
\mbox{\boldmath$\lambda$}_n=\big(an+b+O(n^{-\ep})\big)^{-r}. 
\end{equation*}
\end{Prop}

First, let the boundary conditions be separated or almost separated. As we have proved in Section \ref{Sec2}, see (\ref{separated-final}), the eigenvalues of the first equation in (\ref{equations}) satisfy 
the relation (\ref{asymp}) with $\ep=1$ and $r=3-\alpha$.

To obtain the estimate (\ref{perturb}) we need to normalize $\psi_n$. Since all eigenfunctions except the first one change the sign, the relation (\ref{eigenfun1}) and the estimate (\ref{pogransloy}) imply 
 \begin{equation}\label{eigenfun-est}
|\psi_n(x)|=O(A\nu_n^{-1})=O(An^{-1}), \quad n\to\infty,
\end{equation}
uniformly with respect to $x\in[0,1]$. 

In view of (\ref{norma}) this implies 
\begin{equation*}
\|\psi_n(x)\|_{\cal H}^2=A^2\Big(\int\limits_0^1 \cos^2(\nu_n x)+O(n^{-1})\Big), \quad n\to\infty,
\end{equation*}
so for the normalized eigenfunctions we have $A=\pm\sqrt{2}+O(n^{-1})$. Finally, taking into account (\ref{eigenfun-est}) we obtain
\begin{equation*}
 |({\cal B}\psi_n,\psi_m)_{\cal H}|\le \max\limits_{x\in[0,1]}|\psi_n(x)\psi_m(x)|\int\limits_0^1 |\mathfrak{p}(x)|\,dx \le c(mn)^{-1}
\end{equation*}
for any $\mathfrak{p}\in L_1(0,1)$.\footnote{For some boundary conditions the assumption on $\mathfrak{p}$ can be weakened.
}
Thus, the estimate (\ref{perturb}) is fulfilled with $\ep=1$, and Proposition \ref{prop} ensures the two-term eigenvalues estimate (\ref{separated-final}) for the problem~(\ref{problem1}).

\medskip

For the non-separated boundary conditions the eigenvalues of the operator ${\cal K}$ are organized in two sequences, see (\ref{non-separ-final}). However, the estimate (\ref{perturb}) also holds, and 
Remark 2 in \cite{Naz19} ensures that the asymptotics (\ref{non-separ-final}) persists for the problem (\ref{problem1}).  

Now we can formulate the final result of this section.

\begin{Theorem}\label{pert}
Let $\mathfrak{p}\in L_1(0,1)$. Then the two-term eigenvalue asymptotics of the problem (\ref{problem1}) with self-adjoint boundary conditions does not depend on $\mathfrak{p}$ and is given in 
(\ref{separated-final}) or (\ref{non-separ-final}), depending on boundary conditions.
\end{Theorem}

\section{Gaussian processes related to the problems (\ref{problem}) and (\ref{problem1}) with various boundary conditions}\label{Sec4}

We recall that ${\cal G}$ stands for the covariance function $G_{W^H}$.
\medskip

{\bf 1. Fractional Brownian bridge.} This process is defined as
\begin{equation*}
 B^H(x)=W^H(x)- a(x)W^H(1), \qquad a(x)=\frac {{\cal G}(x,1)}{{\cal G}(1,1)}={\cal G}(x,1). 
\end{equation*}
Its covariance function reads
\begin{equation*}
 G_{B^H}(x,y)={\cal G}(x,y)-{\cal G}(x,1){\cal G}(1,y),
\end{equation*}
and corresponding operator can be considered as a {\it critical} one-dimensional perturbation of the covariance operator of $W^H$, see \cite{Naz09}. 

In \cite{ChiKlM17} this approach was applied to obtain the two-term spectral asymptotics for $B^H$. Moreover, it was mentioned that the direct method developed in \cite{ChiKl} for $W^H$ does not produce 
results quite as explicit as those in \cite{ChiKl}. However, we show that it is not the case, and the direct method works as well.

 Notice that ${G_{B^H}(0,y)=G_{B^H}(1,y)}\equiv0$, and therefore any eigenfunction of
\begin{equation*}
\int\limits_0^1 G_{B^H}(x,y)\varphi(y)\,dy=\lambda\varphi(x)
\end{equation*}
satisfies $\varphi(0)=\varphi(1)=0$.

Define 
\begin{equation}\label{psiBH}
\psi(x)=\int\limits_x^1 \varphi(y)\,dy-c,\qquad c=\int\limits_0^1 {\cal G}(1,y)\varphi(y)\,dy.
\end{equation}
Then evidently $\psi'(0)=\psi'(1)=0$, and
\begin{equation*}
\aligned
\lambda\psi'(x)= &\ \int\limits_0^1 G_{B^H}(x,y)\psi'(y)\,dy\\
= &\ -\int\limits_0^1 ({\cal G}_y(x,y)-{\cal G}(x,1){\cal G}_y(1,y))\psi(y)\,dy.
\endaligned
\end{equation*}
The last term vanishes by the choice of $c$:
\begin{equation}\label{choice-c}
\int\limits_0^1 {\cal G}_y(1,y)\psi(y)\,dy={\cal G}(1,1)\psi(1)+\int\limits_0^1 {\cal G}(1,y)\varphi(y)\,dy=0,
\end{equation}
and we obtain 
\begin{equation*}
\aligned
\lambda\psi'(x)= &\  -\int\limits_0^1 {\cal G}_y(x,y)\psi(y)\,dy\\
= &\ H \int\limits_0^1 \big(y^{2H-1}+\text{sign}(x-y)|x-y|^{2H-1}\big)\psi(y)\,dy.
\endaligned
\end{equation*}
Differentiation gives (\ref{problem}) with $\alpha=2-2H$. Since boundary conditions are separated, we obtain the spectral asymptotics (\ref{separated-final}) with $\varkappa=2$. By Remark \ref{remark-separ}, 
we should exclude zero root of (\ref{eigen-equ}). This changes $n\to n+1$ on the right-hand side of (\ref{separated-final}) and yields
\begin{equation*}
\lambda_n=\sin(\pi H)\Gamma(1+2H) \Big(\pi n+\frac {\pi(H-\frac 12)(\frac 32-H)}{2(H+\frac 12)}+O(n^{-1})\Big)^{-1-2H},
\end{equation*}
that coincides with the result of \cite[Theorem 2.2]{ChiKlM17} and even gives a slightly better estimate of the remainder term.
\medskip

{\bf 2. Centered FBM.} This process is defined as
\begin{equation*}
\overline{W^H}(x)=W^H(x)- \int\limits_0^1 W^H(t)\,dt,  
\end{equation*}
and its covariance function reads
\begin{equation*}
 G_{\overline{W^H}}(x,y)={\cal G}(x,y)-\int\limits_0^1{\cal G}(x,s)\,ds-\int\limits_0^1{\cal G}(t,y)\,dt+\int\limits_0^1\int\limits_0^1{\cal G}(t,s)\,dtds.
\end{equation*}
 Notice that $\int\limits_0^1 G_{\overline{W^H}}(x,y)\,dy=0$. Therefore the equation
\begin{equation*}
\int\limits_0^1 G_{\overline{W^H}}(x,y)\varphi(y)\,dy=\lambda\varphi(x)
\end{equation*}
has a zero eigenvalue corresponding to the constant eigenfunction, and all other eigenfunctions satisfy $\int\limits_0^1\varphi(y)\,dy=0$.

Define $\psi(x)=\int\limits_x^1 \varphi(y)\,dy$.
Then evidently $\psi(0)=\psi(1)=0$, and
\begin{equation*}
\aligned
\lambda\psi'(x)= &\ \int\limits_0^1 G_{\overline{W^H}}(x,y)\psi'(y)\,dy\\
= &\ -\int\limits_0^1 \Big({\cal G}_y(x,y)-\int\limits_0^1{\cal G}_y(t,y)\,dt\Big)\psi(y)\,dy.
\endaligned
\end{equation*}
Differentiation gives (\ref{problem}) with $\alpha=2-2H$. Since boundary conditions are separated, we obtain the spectral asymptotics (\ref{separated-final}) with $\varkappa=0$. This yields
\begin{equation*}
\lambda_n=\sin(\pi H)\Gamma(1+2H) \Big(\pi n-\frac {\pi(H-\frac 12)}2+O(n^{-1})\Big)^{-1-2H}.
\end{equation*}

{\bf 3. Centered Brownian bridge.} This process is defined similarly:
\begin{equation*}
\overline{B^H}(x)=B^H(x)- \int\limits_0^1 B^H(t)\,dt, 
\end{equation*}
and its covariance function reads
\begin{equation*}
\aligned
 G_{\overline{B^H}}(x,y)= &\ G_{B^H}(x,y)-\int\limits_0^1 G_{B^H}(x,s)\,ds\\
 - &\ \int\limits_0^1 G_{B^H}(t,y)\,dt+\int\limits_0^1\int\limits_0^1 G_{B^H}(t,s)\,dtds.
\endaligned
\end{equation*}
 As in the previous example, $\int\limits_0^1 G_{\overline{B^H}}(x,y)\,dy=0$, and thus the equation
\begin{equation}\label{eq-CBH} 
\int\limits_0^1 G_{\overline{B^H}}(x,y)\varphi(y)\,dy=\lambda\varphi(x)
\end{equation}
has a zero eigenvalue corresponding to the constant eigenfunction, while all other eigenfunctions satisfy $\int\limits_0^1\varphi(y)\,dy=0$. Therefore, (\ref{eq-CBH}) can be rewritten as
\begin{equation*}
\int\limits_0^1 \Big(G_{B^H}(x,y)- \int\limits_0^1 G_{B^H}(t,y)\,dt\Big)\varphi(y)\,dy=\lambda\varphi(x),
\end{equation*}
and from ${G_{B^H}(0,y)=G_{B^H}(1,y)}\equiv0$ we conclude that $\varphi(0)=\varphi(1)$.

We define $\psi$ by formula (\ref{psiBH}). Then evidently $\psi(0)=\psi(1)$, $\psi'(0)=\psi'(1)$, and 
\begin{equation*}
\aligned
\lambda\psi'(x)= &\ \int\limits_0^1 \Big(G_{B^H}(x,y)- \int\limits_0^1 G_{B^H}(t,y)\,dt\Big)\psi'(y)\,dy\\
= &\ -\int\limits_0^1 \Big({\cal G}_y(x,y)-{\cal G}(x,1){\cal G}_y(1,y)- \int\limits_0^1 (G_{B^H})_y(t,y)\,dt\Big)\psi(y)\,dy.
\endaligned
\end{equation*}
The second term vanishes by (\ref{choice-c}), and differentiation gives (\ref{problem}) with $\alpha=2-2H$. Since boundary conditions are periodic, we obtain the spectral asymptotics (\ref{non-separ-final}) 
with $\beta=-\gamma$. By Remark \ref{remark-non-separ}, we should exclude zero root of (\ref{eigen-equ}). This changes $n\to n+1$ on the right-hand side of (\ref{non-separ-final}) and yields
\begin{equation*}
\aligned
\lambda_n= &\ \sin(\pi H)\Gamma(1+2H) 
\Big(\pi n-\frac {\pi(H-\frac 12)}2+\frac {\pi}2\,(1-(-1)^n)\\
- &\ \frac {\pi(H-\frac 12)}{2(H+\frac 12)}\,(1-(-1)^n)+O(n^{-1})\Big)^{-1-2H}.
\endaligned
\end{equation*}

{\bf 4. Fractional Slepian process.} The conventional Slepian process $S$ on $[0,1]$ can be defined in several ways:

1) $S$ is a zero mean-value stationary Gaussian process with correlation function $1-|x-y|$;

2) $S(x)=W(x+1)-W(x)$ where $W$ is the Wiener process;

3) $S(x)=W_1(x)+W_2(1-x)$ where $W_1$ and $W_2$ are independent Wiener processes.

Fractional Slepian processes defined by analogy in these three ways are different, see, e.g. \cite{Mol} for the first one. We define the process $S^H$ as the mixture of two independent FBMs:
\begin{equation*}
 S^H(x)=W_1^H(x)+W_2^H(1-x),\qquad x\in[0,1].
\end{equation*}
 Its covariance function reads
\begin{equation*}
 G_{S^H}(x,y)= {\cal G}(x,y)+{\cal G}(1-x,1-y).
\end{equation*}
We emphasize that, as in the case $H=1/2$, the following relation holds:
\begin{equation*}
 S^H(x)-S^H(0)\stackrel{d}{=}2W^H(x).
\end{equation*}
Notice also that $G_{S^H}(0,y)+G_{S^H}(1,y)\equiv1$.

By symmetry of the kernel, any eigenfunction of
\begin{equation*}
\int\limits_0^1 G_{S^H}(x,y)\varphi(y)\,dy=\lambda\varphi(x)
\end{equation*}
satisfies either $\varphi(x)=\varphi(1-x)$ or $\varphi(x)=-\varphi(1-x)$.

In the first case we define 
\begin{equation}\label{psiSH1}
\psi(x)=\int\limits_{\frac 12}^x \varphi(y)\,dy.
\end{equation}
Evidently $\psi(0)+\psi(1)=0$, and 
\begin{equation*}
\lambda(\varphi(0)+\varphi(1))=\int\limits_0^1 \big(G_{S^H}(0,y)+G_{S^H}(1,y)\big)\varphi(y)\,dy=
\int\limits_0^1 \varphi(y)\,dy,
\end{equation*}
that is equivalent to
\begin{equation}\label{bcSH}
\psi(0)+\psi(1)=0;\qquad \psi'(0)+\psi'(1)+\frac 2{\lambda}\psi(0)=0.
\end{equation}

In the second case we define
\begin{equation}\label{psiSH2}
\psi(x)=\int\limits_0^x \varphi(y)\,dy.
\end{equation}
Then evidently $\psi(0)=\psi(1)=0$ and $\psi'(0)+\psi'(1)=0$, so the boundary conditions (\ref{bcSH}) are satisfied as well.

Further,
\begin{equation*}
\aligned
\lambda\psi'(x)= &\ \int\limits_0^1 G_{S^H}(x,y)\psi'(y)\,dy=G_{S^H}(x,1)\psi(1)-G_{S^H}(x,0)\psi(0)\\
- &\ \int\limits_0^1 ({\cal G}_y(x,y)-{\cal G}_y(1-x,1-y))\psi(y)\,dy.
\endaligned
\end{equation*}
The double substitution is equal to $\psi(1)$ and vanishes after differentiation. So, we obtain 
\begin{equation}\label{problem2}
(\mathbb{K}_{\alpha}\psi)(x)=-\frac {\lambda}2\, \psi''(x)
\end{equation}
with $\alpha=2-2H$ and boundary conditions (\ref{bcSH}).

The boundary conditions in this problem are non-separated but contain the spectral parameter $\lambda$. So, formula (\ref{non-separ-final}) is not applicable. However, the basic scheme runs without 
essential changes. We change $\lambda\to \lambda/2$ in (\ref{nu}) and arrive at (\ref{eigen3}). Using (\ref{bcSH}) we rewrite (\ref{eigen3}) as follows:
\begin{equation*}
\aligned
& 2 \big[\mathfrak{A}\big]\psi'(0)+\nu^{3-\alpha}\frac{\cos(\pi\alpha/2)}{c_\alpha\pi}\big[\mathfrak{A}\big]\psi(0)=0;\\
& \Big(\nu^{3-\alpha}\frac{\cos(\pi\alpha/2)}{c_\alpha\pi}\big[\mathfrak{B}\big]-2\nu\big(\big[\mathfrak{A}\big]+\mathfrak{b}_{\alpha}\big[\mathfrak{B}\big]\big)\Big)
\psi(0)=0.
\endaligned
\end{equation*}
The existence condition for a nontrivial solution to this system takes the form
\begin{equation*}
2\mathfrak{A}\mathfrak{B}\equiv \sin(\nu+\rho)=O(\nu^{-\min\{1,2-\alpha\}}),\qquad \mbox{as}\ \ \nu\to\infty,
\end{equation*}
and we conclude
\begin{equation*}
\nu_{n+k}=\pi n-\frac {\pi(1-\alpha)}4+O(n^{-\min\{1,2-\alpha\}}),\qquad \mbox{as}\ \ n\to\infty
\end{equation*}
for some $k$. Comparing this result with the eigenvalues of the conventional Slepian process, see \cite{NiO}, we obtain $k=1$, and
\begin{equation*}
\lambda_n= 2\sin(\pi H)\Gamma(1+2H)\Big(\pi (n-1)-\frac {\pi(H-\frac 12)}2+O(n^{-\min\{1,2H\}})\Big)^{-1-2H}\!.
\end{equation*}

{\bf 5.} The spectral analysis of some other mixtures of fractional processes can be reduced to the problems of the same type, cf. \cite[Section 2]{NaNi18} for the mixtures of the Green Gaussian processes. 
We consider here only one modification of the previous example.

Let $S^H_\gamma(x)=S^H(x)-\gamma(S^H(0)+S^H(1))$. This process is a one-dimensional perturbation of the fractional Slepian process, and its covariance function reads
\begin{equation*}
 G_{S^H_\gamma}(x,y)= {\cal G}(x,y)+{\cal G}(1-x,1-y)+2(\gamma^2-\gamma).
\end{equation*}
For $\gamma\ne\frac 12$ this perturbation is {\it non-critical}, see \cite{Naz09}, and the two-term spectral asymptotics does not change. The case $\gamma=\frac 12$ is critical and should be studied separately, 
cf. \cite[Theorem 2.1]{Naz09a} and \cite[Example 8]{Naz09} for conventional Slepian process.

To manage the equation for eigenfunctions in the case $\gamma=\frac 12$
\begin{equation*}
\int\limits_0^1 \big(G_{S^H}(x,y)-\frac 12\big)\varphi(y)\,dy=\lambda\varphi(x)
\end{equation*}
we notice that $\varphi(0)+\varphi(1)=0$ and make the change of function (\ref{psiSH1}), (\ref{psiSH2}). This gives (\ref{problem2}) with $\alpha=2-2H$ and anti-periodic boundary conditions
\begin{equation*}
\psi(0)+\psi(1)=0;\qquad \psi'(0)+\psi'(1)=0.
\end{equation*}
Thus we obtain the spectral asymptotics (\ref{non-separ-final}) with $\beta=\gamma$. This yields
\begin{equation*}
\aligned
\lambda_n= &\ 2\sin(\pi H)\Gamma(1+2H) 
\Big(\pi n-\frac {\pi(H-\frac 12)}2-\frac {\pi}2\,(1+(-1)^n)\\
+ &\ \frac {\pi(H-\frac 12)}{2(H+\frac 12)}\,(1+(-1)^n)+O(n^{-1})\Big)^{-1-2H}.
\endaligned
\end{equation*}

{\bf 6. The fractional Ornstein--Uhlenbeck process beginning at zero.} In the fractional setting, the Ornstein--Uhlenbeck process can be defined in a number of nonequivalent ways, see, e.g., \cite{CheKaM}. 
Following \cite{ChiKlM17}, we consider the solution of the Langevin equation driven by the FBM:
\begin{equation}\label{OU-def}
 U_\beta^H(x)=\xi-\beta\int\limits_0^x U_\beta^H(t)\,dt+W^H(x),
\end{equation}
where $\beta\in\mathbb{R}$ is the drift parameter and $\xi\sim {\cal N}(0,\sigma^2)$ is the initial condition independent of $W^H$. The covariance function is given by the formula
\begin{equation*}
\aligned
 G_\beta(x,y)\equiv &\ G_{U_\beta^H}(x,y)=\exp(-\beta(x+y))\\
 \times\Big[\sigma^2+ & \int\limits_0^x\exp(\beta s)\frac d{ds}\int\limits_0^y H|s-t|^{2H-1}\text{sign}(s-t)\exp(\beta t)\,dtds\Big].
\endaligned
\end{equation*}

We begin with $\sigma=0$, i.e. $\xi=0$. This case was considered in \cite[Sec. 6]{ChiKlM18}. By certain fine analysis the following expression for eigenvalues was derived:
\begin{equation}\label{OU-ChKlM}
\lambda_n=\sin(\pi H)\Gamma(1+2H)\, \frac {\nu_n^{1-2H}}{\nu_n^2+\beta^2},\qquad n\in\mathbb{N},
\end{equation}
where the sequence $\nu_n$ satisfies (\ref{kappa=1}) with $k=0$ and $\alpha=2-2H$.

We claim that this result is covered by our Theorem \ref{pert}. Indeed, the change of function 
\begin{equation}\label{OU-change}
\psi(x)=\exp(\beta x)\int\limits_x^1 \exp(-\beta y)\varphi(y)\,dy
\end{equation}
(cf.~the proof of Lemma 6.1 in \cite{ChiKlM18}) reduces the equation
\begin{equation*}
\int\limits_0^1 G_{\beta}(x,y)\varphi(y)\,dy=\lambda\varphi(x)
\end{equation*}
to the problem
\begin{equation}\label{OU-probl}
(\mathbb{K}_{\alpha}\psi)(x)=\lambda \big(-\psi''(x)+\beta^2\psi(x)\big),\quad x\in(0,1)
\end{equation}
with boundary conditions $(\psi'-\beta\psi)(0)=\psi(1)=0$.

Theorem \ref{pert} shows that 
\begin{equation*}
\lambda_n= \sin(\pi H)\Gamma(1+2H) 
\Big(\pi n-\frac {\pi(H-\frac 12)}2-
\frac {\pi}{2(H+\frac 12)}+O(n^{-1})\Big)^{-1-2H},
\end{equation*}
which coincides with (\ref{OU-ChKlM}) taking into account that $\nu_n^2+\beta^2=\nu_n^2(1+O(n^{-2}))$. Thus, the claim follows.
\medskip

{\bf 7.} Now we consider the Ornstein--Uhlenbeck process (\ref{OU-def}) with $\sigma\ne0$. The change of function (\ref{OU-change}) gives the same equation (\ref{OU-probl}) and the boundary condition $\psi(1)=0$. 
To obtain the boundary condition at zero we write
\begin{equation*}
\lambda\varphi(0)=\int\limits_0^1 G_{\beta}(0,y)\varphi(y)\,dy=\sigma^2 \int\limits_0^1 \exp(-\beta y)\varphi(y)\,dy=\sigma^2\psi(0),
\end{equation*}
and the evident relation $\beta\psi(x)-\psi'(x)=\varphi(x)$ implies
\begin{equation}\label{bcOU}
\psi'(0)-\big(\beta-\frac {\sigma^2}{\lambda}\big)\psi(0)=0;\qquad \psi(1)=0.
\end{equation}

First, we consider the problem (\ref{problem}) with the same boundary conditions (\ref{bcOU}). These boundary conditions are separated but contain the spectral parameter $\lambda$. So, as in the example {\bf 4}, 
formula (\ref{separated-final}) is not applicable. However, the basic scheme again runs without changes. Using (\ref{bcOU}) we rewrite (\ref{eigen3}) as follows:
\begin{equation*}
\aligned
& \Big(\sigma^2\nu^{3-\alpha}\frac{\cos(\pi\alpha/2)}{c_\alpha\pi}\big[\mathfrak{A}\big]-\nu\big(\mathfrak{b}_{\alpha}\big[\mathfrak{A}\big]-\big[\mathfrak{B}\big]\big)-\beta\big[\mathfrak{A}\big]\Big)\psi(0)
+ \big[\mathfrak{A}\big]\psi'(1)=0;\\
& \Big(\sigma^2\nu^{3-\alpha}\frac{\cos(\pi\alpha/2)}{c_\alpha\pi}\big[\mathfrak{B}\big]-\nu\big(\big[\mathfrak{A}\big]+\mathfrak{b}_{\alpha}\big[\mathfrak{B}\big]\big)-\beta\big[\mathfrak{B}\big]\Big)
\psi(0)-\big[\mathfrak{B}\big]\psi'(1)=0.
\endaligned
\end{equation*}
The existence condition for a nontrivial solution to this system takes the form
\begin{equation*}
2\mathfrak{A}\mathfrak{B}\equiv \sin(\nu+\rho)=O(\nu^{-\min\{1,2-\alpha\}}),\qquad \mbox{as}\ \ \nu\to\infty,
\end{equation*}
and we conclude
\begin{equation*}
\nu_{n+k}=\pi n-\frac {\pi(1-\alpha)}4+O(n^{-\min\{1,2-\alpha\}}),\qquad \mbox{as}\ \ n\to\infty
\end{equation*}
for some $k$. 

We postpone the specification of $k$ and turn to the problem (\ref{OU-probl}). Theorem \ref{pert} is not applicable directly because of the spectral parameter in the boundary conditions (\ref{bcOU}). However, 
basic relations (\ref{eigenfun1}) and (\ref{pogransloy}), and therefore the estimate (\ref{eigenfun-est}) hold regardless of the boundary conditions. To apply Proposition \ref{prop} we need only to redefine the operator ${\cal K}$ in (\ref{equations}) denoting
\begin{equation*}
({\cal K}\psi,\eta)_{\cal H}:=\int\limits_0^1(\mathbb{K}_{\alpha}\psi)(x)\eta(x)\,dx +\sigma^2\psi(0)\eta(0).
\end{equation*}
The estimate (\ref{perturb}) with $\ep=1$ persists, and Proposition \ref{prop} shows that the term $\beta^2u$ in (\ref{OU-probl}) does not affect the two-term eigenvalue asymptotics. 

Finally, comparing this result with the eigenvalues of the conventional Ornstein--Uhlenbeck process corresponding to $\alpha=1$ and $\sigma^2=\frac 1{2\beta}$, see, e.g., \cite[Proposition 5.5]{NaNi04}, 
we obtain $k=1$, and
\begin{equation*}
\lambda_n= \sin(\pi H)\Gamma(1+2H)\Big(\pi (n-1)-\frac {\pi(H-\frac 12)}2+O(n^{-\min\{1,2H\}})\Big)^{-1-2H}\!.
\end{equation*}

\section{Application to the small ball probabilities}\label{Sec5}

As it was mentioned in the Introduction, the results of the previous section give rise to the exact (up to a constant) $L_2$-small ball asymptotics of all considered Gaussian processes. 

We define two important quantities:
\begin{equation*}
\aligned
D(H):= &\ \frac H{(2H+1)\sin\left(\frac {\pi}{2H+1}\right)} \bigg (\frac {\sin(\pi H)\Gamma(2H+1)} {(2H+1)\sin\left(\frac {\pi}{2H+1}\right)}\bigg)^{\frac {1}{2H}};\\
B(H):= &\ \frac {(H-\frac 12)^2}{2H}.
\endaligned
\end{equation*}

We substitute the two-term spectral asymptotics into the general result of \cite[Theorem 6.2]{NaNi04}. For the examples {\bf 3} and {\bf 5} we additionally use the Lifshits lemma, see, e.g., 
\cite[Lemma 0.1]{Naz09a}. This gives us the following statement.

\begin{Theorem}
The exact small ball asymptotics for the fractional processes considered in Section \ref{Sec4} read as follows:
\begin{equation*}
\mathbb{P}\big\{\int\limits_0^1 X^2(x)\,dx\leq\ep^2\big\} \sim C(X)\cdot \ep^{B_X} \exp(- D_X \ep^{-\frac 1H}), 
\end{equation*}
where the values of $B_X$ and $D_X$ are collected in the table.\footnote{The result for $W^H$ was derived in \cite{ChiKl} and is given for the comparison.}
\begin{center}
\begin{tabular}{|c|c|c|} 
\hline
$X$ & $B_X$ & $D_X\vphantom{\frac {1^1}2}$ \\
\hline \hline 
$W^H$ & $B(H)+\frac 1{2H}$ & $D(H)\vphantom{\frac {1^{1^1}}{2^{2^2}}}$\\ 
\hline 
$B^H$ & $B(H)+\frac 1{2H}-1$ & $D(H)\vphantom{\frac {1^{1^1}}{2^{2^2}}}$\\ 
\hline 
$\overline{W^H}$ & $B(H)$ & $D(H)\vphantom{\frac {1^{1^1}}{2^{2^2}}}$\\ 
\hline
$\overline{B^H}$ & $B(H)-1$ & $D(H)\vphantom{\frac {1^{1^1}}{2^{2^2}}}$\\ 
\hline 
$S^H_{\gamma}, \gamma\ne\frac 12$ & $B(H)+\frac 1{2H}+1$ & $2^{\frac 1{2H}}D(H)\vphantom{\frac {1^{1^1}}{2^{2^2}}}$\\ 
\hline
$S^H_{\frac 12}$ & $B(H)+\frac 1{2H}$ & $2^{\frac 1{2H}}D(H)\vphantom{\frac {1^{1^1}}{2^{2^2}}}$\\ 
\hline 
$U^H_{\beta}$, $\sigma=0$ & $B(H)+\frac 1{2H}$ & $D(H)\vphantom{\frac {1^{1^1}}{2^{2^2}}}$\\ 
\hline
$U^H_{\beta}$, $\sigma\ne0$ & $B(H)+\frac 1{2H}+1$ & $D(H)\vphantom{\frac {1^{1^1}}{2^{2^2}}}$\\ 
\hline
\end{tabular}
\end{center}
\smallskip
\end{Theorem}

\begin{Remark}
1. It is well known that the centered Wiener process coincides in distribution with the Brownian bridge. The table shows that for $H\ne\frac 12$ this is not the case, and even $L_2$-small ball asymptotics 
for $B^H$ and $\overline{W^H}$ differ at the power level.

2. In contrast, $L_2$-small ball asymptotics for $W^H$ and $U^H_{\beta}$ in the case $\sigma=0$ coincide up to a constant for all $H\in(0,1)$.
\end{Remark}

\paragraph*{Acknowledgements.}
I am greatful to P.~Chigansky, M.~Kleptsyna and Ya.~Nikitin for useful discussions. Also I would like to thank N.~Rastegaev and the anonymous referee for their careful reading and some important advice.

The author's work is supported by the joint RFBR--DFG grant 20-51-12004.

\end{document}